\documentclass[letterpaper, 10 pt, conference]{ieeeconf}  
\IEEEoverridecommandlockouts                              
\usepackage{amsfonts, amsmath, bm}
\usepackage{amssymb}
\usepackage{algorithmic, graphicx}
\usepackage{epstopdf}
\usepackage{float}
\usepackage{graphicx}
\usepackage{hyperref}
\usepackage{tikz}
\usepackage{pgfplots}
\pgfplotsset{compat=1.18}
\usetikzlibrary{calc}
\usetikzlibrary{arrows.meta}
\usepackage{subcaption}
\usepackage{color}
\usepackage{xcolor}

\title{\LARGE \bf
Extended State Observer for Localized Fault Awareness in RF Accelerating Structures
} 
\author{Luke S. Baker, Sungil Kwon, Kwame Jyamfi, Quinten Cole, Isaac Roybal, AJ Garcia, Phil Torrez,\\ Lawrence Castellano
\thanks{The authors would like to thank Alexander Scheinker, Anatoly Zlotnik, and Charles E. Taylor for helpful discussions.  All authors are with the Accelerator Operations and Technology Division at the Los Alamos National Laboratory, Los Alamos, NM, USA 87545.  L.S. Baker (lsbaker@lanl.gov) is the author of correspondence. Research conducted at Los Alamos National Laboratory is done under the auspices of the National Nuclear Security Administration of the U.S. Department of Energy under Contract No. 89233218CNA000001. Report No. LA-UR-26-22463.}}

\begin{document}

\maketitle
\thispagestyle{empty}
\pagestyle{empty}

\begin{abstract}
An observer framework is presented for robust regulation of RF cavity fields and localized identification of disturbances in RF systems. A standard cavity field observer is augmented with additional states to estimate the evolution of cavity detuning and phase drifts induced by the drive and receiver chains. Monte Carlo simulations are performed to assess the performance of the proposed estimator under realistic conditions for the intended high-power linear accelerator operation. Results showcase precise cavity field regulation and the reliability with which the observer assigns deviations to the correct subsystem. The resulting diagnostic capability provides a foundation for improved fault detection, faster troubleshooting during accelerator operation, and more informed maintenance of RF systems in large accelerator facilities.
\end{abstract}

\section{Introduction}
\label{sec:introduction}

Radio-frequency (RF) accelerating cavities are fundamental components for transferring energy in modern particle accelerators.  These resonant electromagnetic structures store RF energy synchronized with the arrival of charged particle bunches to efficiently transfer energy from high-power RF amplifiers such as klystrons and solid-state transmitters to the beam \cite{wangler2008,padamsee2009}. Chains of coupled resonant cavities form accelerating structures through which the electromagnetic behavior is accurately modeled using the theory of coupled resonators \cite{nagel1964coupled,nagle1967}. These models describe dispersion characteristics, field stability, and sensitivity to fabrication errors in multicell cavity structures. Because the energy gain imparted to the beam depends directly on the amplitude and phase of the accelerating field, maintaining stable cavity fields is essential for preserving quality beam and minimizing energy spread during accelerator operation \cite{simrock2022low}.

\begin{figure}
\centering
\begin{tikzpicture}[
    >=stealth,
    thick,
    x=0.78cm,
    y=0.7cm,
    bunch/.style={
        fill=blue!70!black,
        draw=blue!70!black,
        rounded corners=0.6pt
    },
    envline/.style={red!80!black, ultra thick},
    sample/.style={red!80!black, fill=red!80!black}
]

\def\xFillEnd{2.8}      
\def\xFlatEnd{7.8}      
\def\xEnd{9.5}          

\draw[->,blue!70!black] (0,0) -- (0,3.1)
    node[pos=1.0, anchor=center, left=12pt, rotate=90, text=blue!70!black]
    {Beam Current};

\draw[->,red!80!black] (\xEnd,0) -- (\xEnd,3.1)
    node[pos=0.0, anchor=center, right=12pt, rotate=90, text=red!80!black]
    {RF Envelope};

\draw[->,black!85] (0,0) -- (\xEnd+0.25,0);
\draw[blue!70!black, line width=1.2pt] (0,0) -- (\xEnd,0);

\foreach \x/\lab in {0/0,\xFillEnd/325~\mu s,\xFlatEnd/950~\mu s}{
  \draw[black!80] (\x,0.08) -- (\x,-0.08);
  \node[black!80, below] at (\x,-0.08) {\scriptsize $\lab$};
}

\draw[gray!55, dashed] (\xFillEnd,0) -- (\xFillEnd,3.0);
\draw[gray!55, dashed] (\xFlatEnd,0) -- (\xFlatEnd,3.0);

\node[gray!65] at ({0.5*\xFillEnd},2.95) {\scriptsize Fill};
\node[gray!65] at ({0.5*(\xFillEnd+\xFlatEnd)},2.95) {\scriptsize Flattop};
\node[gray!65] at ({0.5*(\xFlatEnd+\xEnd)},2.95) {\scriptsize RF Off};

\draw[<->,gray!70] (0,-0.55) -- (\xFillEnd,-0.55)
    node[midway,below=2pt]{ $325~\mu s$};

\draw[<->,gray!70] (\xFillEnd,-0.55) -- (\xFlatEnd,-0.55)
    node[midway,below=2pt]{ $625~\mu s$};

\draw[<->,gray!70] (\xFlatEnd,-0.55) -- (\xEnd,-0.55)
    node[midway,below=2pt]{ $7.38~ms$};

\draw[gray!70,line width=0.5pt] ({\xFlatEnd+0.55},-0.64) -- ({\xFlatEnd+0.70},-0.46);
\draw[gray!70,line width=0.5pt] ({\xFlatEnd+0.78},-0.64) -- ({\xFlatEnd+0.93},-0.46);

\def\jitter{0.15}

\foreach \x in {3.2,3.6,...,7.7}{
    \pgfmathsetmacro{\h}{2.15 + \jitter*(2*rnd - 1)}
    \draw[bunch] (\x-0.025,0) rectangle (\x+0.025,\h);
}
\foreach \x in {3.1,3.2,...,7.7}{
    \pgfmathsetmacro{\h}{2.15 + \jitter*(2*rnd - 1)}
    \draw[bunch] (\x-0.0008,0) rectangle (\x+0.0008,\h);
}


\draw[dashed,red!80!black] (0,1.8) -- ({\xEnd-0.90},1.8);
\node[red!80!black] at ({\xEnd-0.50},1.8) {$V^{\mathrm{cav}}$};

\draw[envline,domain=0:\xFillEnd,smooth,variable=\x]
    plot ({\x},{1.8*(1 - exp(-5*\x/\xFillEnd))});

\draw[envline,domain=\xFillEnd:\xFlatEnd,smooth,variable=\x]
    plot ({\x},{1.8*(1 - exp(-5*\x/\xFillEnd)) + 0.06*sin(1.9*(\x -\xFillEnd) r)});

\draw[envline,domain=\xFlatEnd:\xEnd,smooth,variable=\x]
    plot ({\x},{1.8*exp(-7.0*(\x-\xFlatEnd))});

\end{tikzpicture}
\caption{Illustration of a macropulse repeated at 120 Hz. Beam is present only during the regulated flattop with a bunch repetition rate of 805 MHz for the side-coupled cavities at LANSCE.}
\label{fig:macropulse-single}
\end{figure}
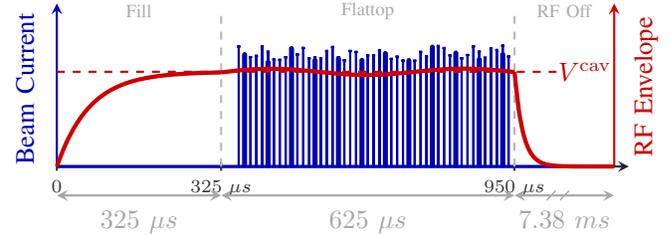

The Los Alamos Neutron Science Center (LANSCE) provides a prominent example of an RF linear accelerator (linac) operating under such conditions. The LANSCE linac accelerates protons to energies approaching 800~MeV using a sequence of RF accelerating structures operating at multiple frequencies \cite{swenson1973alamos,young2007lansce}. Following particle injection and initial bunching, acceleration begins with a 201.25 MHz drift-tube linac (DTL) that raises beam energy to approximately 100 MeV \cite{young2007lansce}.  This is followed by a series of 805 MHz side-coupled cavity modules that provide the remaining energy gain \cite{swain1973cavity}. The side-coupled architecture enables stable $\pi/2$-mode operation in long cavity chains while maintaining efficient RF coupling and field uniformity across the structure \cite{swain1973cavity,swenson1967stabilization}. Over several decades of operation, LANSCE has provided extensive experience with the long-term operational behavior of RF systems in demanding accelerator environments \cite{jameson1975lampf,lyles2016design}.  As illustrated in Figure \ref{fig:macropulse-single}, LANSCE operates in pulsed mode up to 120 Hz, where RF power is enabled only during the fill and flattop portions of each macropulse.

A distinctive feature of LANSCE operation is the wide diversity of beam flavors delivered to multiple experimental facilities, including the Lujan Neutron Scattering Center, the Proton Radiography Facility, the Weapons Neutron Research Facility, the Isotope Production Facility, and the Ultra Cold Neutron Facility \cite{scheinker2017iterative,taylor2025post}. These programs require beam pulses with substantially different average currents, repetition rates, and time structures \cite{alvinerie2025lansce}. As a result, the RF cavities experience a wide range of beam-loading conditions and disturbances that depend on beam current and its phase relative to the RF waveform \cite{boussard2007control,scheinker2021extremum}. Additional perturbations arise from cavity detuning caused by thermal expansion, Lorentz-force deformation, and mechanical vibrations.  These shift the resonant frequency and modify the efficiency of RF power transfer to the cavity fields \cite{leewe2017resonance}.  In addition to intrinsic cavity perturbations, temperature variations in transmission cables and electronic components can shift the phases of RF drive and pickup signals relative to the reference.  Projects at LANSCE continue to upgrade RF infrastructure and digital low-level RF (LLRF) control systems to improve stability and meet the demands of an increasingly wide range of beam conditions \cite{rees2007lansce,prokop2007lansce,van2022lansce,lyles2022rf}.

Early accelerator facilities employed analog feedback loops and feedforward compensation to mitigate beam-loading transients and maintain acceptable field stability in pulsed linacs. Advances in digital signal processing and field-programmable gate array (FPGA) platforms have since enabled more flexible and precise digital LLRF architectures \cite{giergusiewicz2005fpga,kwon2016fpga,wibowo2018digital}. Modern systems represent field signals in complex baseband (IQ) form derived by approximate RLC circuit models to enable vector control of cavity fields using high-speed digital feedback loops \cite{adaptiveRF2002,scheinker2017preliminary,kwon2011fpga,jablonski20152pi}. These architectures integrate RF detection, digital control algorithms, and timing synchronization to achieve amplitude and phase stability with high precision. 
Large accelerator projects such as the International Linear Collider and the European X-Ray Free Electron Laser (XFEL) have more recently provided opportunity to drive advances in RF technology for high-precision RF reference distribution, advanced RF diagnostics, and energy upgrades beyond 1 to 3 TeV \cite{schilcher1998vector,simrock2005low,vogel2007high,shemelin2022optimization}.

Proportional–integral (PI) control combined with feedforward compensation is widely used in digital LLRF systems as the basis for feedback regulation of cavity field amplitude and phase \cite{champion2003spallation,simrock2022low,kwon2011fpga}. Robust and model-predictive control techniques have been proposed to reject disturbances and incorporate cavity dynamics more directly into the control design \cite{pfeiffer2012design,xiangping2014model,wang2023enhancing}. At the same time, adaptive and model-independent control methods have been explored to address time-varying operating conditions and uncertainties in RF cavity systems. Techniques based on extremum seeking have been applied to automatically tune cavity resonance and optimize RF system performance without requiring detailed models of the underlying dynamics \cite{7526624,scheinker2016application}. Adaptive resonance control and beam-loading compensation strategies have also been proposed for operation under changing beam conditions \cite{scheinker2017iterative,scheinker2021extremum}. These developments reflect a broader trend toward RF control systems capable of accommodating multiple sources of variation while providing improved insight into system behavior \cite{scheinker2021adaptive,simrock2022low}.

A persistent challenge in RF cavity control is distinguishing between physical sources of phase and amplitude error that produce similar signatures in measured RF signals. In digital LLRF systems, cavity detuning and phase drifts in the forward drive and receiver chains can all appear as rotations of the cavity phasor in the complex baseband representation \cite{liepe2001dynamic,yamamoto2009observation,richter2025estimation,doolittle2023drift}. Conventional feedback controllers regulate the cavity field without explicitly identifying the origin of these disturbances, which complicates diagnostics and fault detection during accelerator operation. This work develops an observer framework that separates these disturbance mechanisms within a unified baseband cavity model by augmenting the standard cavity field observer with additional states that estimate forward-chain phase drift, cavity detuning, and receiver-chain phase drift. Beyond cavity field regulation, this approach provides a means to identify both the source and location of disturbances within the RF system for more informed diagnostics and troubleshooting before performance degrades or protective interlocks are triggered.

The remainder of the paper is organized as follows. Section~\ref{sec:baseband-model} develops a complex baseband model of the cavity dynamics, RF drive chain, and receiver system. Section~\ref{sec:estimation} introduces the observer architecture used to estimate cavity state, additive disturbances, detuning, and hardware phase drifts. Section~\ref{sec:control-architecture} presents the feedforward and feedback control strategies used to regulate the cavity field based on these estimates. Section~\ref{sec:mc_fault_awareness} evaluates the performance of the proposed approach using Monte Carlo simulations that compare regulation performance and fault localization capability against a conventional observer structure.

\section{Baseband Cavity Model}
\label{sec:baseband-model}

\begin{figure}[b]
\centering
\resizebox{\linewidth}{!}{%
\begin{tikzpicture}[x=1mm,y=1mm, font=\small, >=Latex, thick]
\tikzset{
  box/.style={draw, rounded corners, align=center, minimum height=7mm, fill=blue!4},
  proc/.style={draw, rounded corners, align=center, minimum height=7mm, fill=green!6},
  dist/.style={draw, rounded corners, align=center, minimum height=6mm, fill=orange!12},
  sum/.style={draw, circle, inner sep=1.1pt, fill=white},
  line/.style={->, thick},
  dline/.style={->, thick, dashed}
}

\node[box, minimum width=16mm] (ref)  at (20,16) {$x_k^*$};
\node[sum] (s1) at (20,0) {$\pm$};

\node[proc, minimum width=20mm] (ctl)  at (40,0) {FCM};
\node[box,  minimum width=20mm] (drv)  at (71,0) {Forward \\ chain};
\node[dist, minimum width=18mm] (fwd) at (71,16) {$\varphi^{\mathrm{fwd}}_k$};
\node[proc, minimum width=26mm] (cav)  at (102,0) {Cavity\\$A_k,B_k$};
\node[dist, minimum width=18mm] (rec) at (102,-32) {$\varphi^{\mathrm{rec}}_k$};

\node[box, minimum width=22mm] (meas) at (102,-18) {Receiver \\ chain};
\node[sum] (sy) at (77,-18) {$+$};
\node[dist, minimum width=18mm] (noi) at (77,-32) {$\nu_k$};

\node[sum] (sr) at (60.6,-18) {$\pm$};
\node[proc, minimum width=22mm] (obs) at (40,-18) {ESO};

\node[proc, minimum width=22mm] (pred) at (40,-32) {$\hat x_k^{\rm pred}$};

\draw[line] (ref) -- (s1);
\draw[line] (s1) -- node[pos=0.4,above] {$\hat e_k$} (ctl);
\draw[line] (ctl) -- node[pos=0.4,above] {$u_k$} (drv);
\draw[line] (drv) -- node[above] {$u_k^{\mathrm{fwd}}$} (cav);
\draw[dline] (fwd.south) -- (drv.north);
\draw[dline] (rec.north) -- (meas.south);

\draw[line] (cav) -- node[right] {$x_k$} (meas);

\draw[line] (meas.west) -- (sy.east);
\draw[dline] (noi.north) -- (sy.south);

\draw[line] (sy.west) -- node[pos=0.4,below] {$y_k$} (sr.east);

\draw[line] (pred.east) -| (sr.south);

\draw[line] (sr.west) -- node[pos=0.25,below] {$r_k$} (obs.east);

\draw[line] (obs.west) -| node[pos=0.25,below] {$\hat x_k$ } (s1.south);

\node[dist, minimum width=18mm] (det)  at (102,16) {$\Delta\omega_k,\, d_k$};
\draw[dline] (det.south) -- (cav.north);
\draw[dline] (obs.north) -- node[right] {$\hat d_k$} (ctl.south);

\end{tikzpicture}%
}
\caption{Simplified block diagram of the LLRF control system at LANSCE.  Orange input boxes represent sources of disturbance, uncertainty, and noise.}
\label{fig:arch-highlevel-compact}
\end{figure}
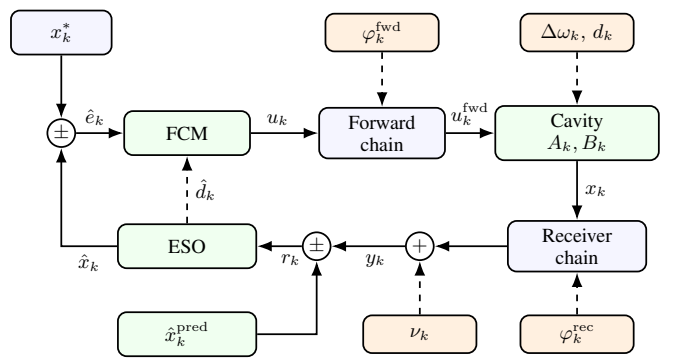

An accelerating cavity operating in its fundamental mode is
approximated as a narrowband resonator at nominal angular frequency
$\omega_0$ and loaded quality factor $Q_L$.
To model the cavity dynamics and associated control system,
all voltages are represented as complex baseband envelopes in a frame
rotating at the LLRF reference angular frequency
$\omega^{\mathrm{ref}}$ \cite{piller2005spallation}.
The physical accelerating voltage is given by
\begin{equation}
  v^{\mathrm{cav}}(t) = \Re \left\{ V^{\mathrm{cav}}(t)\,e^{j\omega^{\mathrm{ref}} t} \right\},
\end{equation}
where $V^{\mathrm{cav}}(t)$ is a slowly varying complex envelope whose magnitude
and phase determine the accelerating gradient and RF phase relative
to the LLRF reference. This section derives a complex baseband state-space model for
estimating $V^{\mathrm{cav}}(t)$. The objective is to separate slow drifts
introduced by the forward drive and receiver measurement chains from detuning
dynamics associated with cavity operation off resonance.  Figure \ref{fig:arch-highlevel-compact} illustrates the LLRF feedback loop and channels through which disturbances enter the system.

\subsection{Drive Synthesis and Forward Path}

In digital LLRF systems such as those at LANSCE, a complex baseband voltage command $V(t)$ is synthesized by an FPGA within a field control module (FCM) or equivalent hardware. This signal is converted by a DAC (digital-to-analog converter) to an analog intermediate frequency (IF), 
upconverted to the RF carrier, amplified by the high-power amplifier chain, 
and delivered to the cavity input coupler by waveguide or coaxial transmission. Temperature variations and slow phase drift in the modulator, upconverter, 
amplifiers, and transmission lines introduce a phase rotation between the 
commanded signal and the RF drive that reaches the cavity. We model this hardware 
drift in the forward drive chain as a lumped phase 
rotation $\varphi^{\mathrm{fwd}}(t)$ between the commanded drive $V(t)$ and 
the forward voltage $V^{\mathrm{fwd}}(t)$ incident on the cavity according to
\begin{equation}
V^{\mathrm{fwd}}(t) = e^{j\varphi^{\mathrm{fwd}}(t)}\,V(t).
\end{equation}
Although the forward coupler may not be exactly incident on the cavity, any additional phase drift accumulated along the 
transmission segment between the coupler and the cavity is considered negligible 
compared to that introduced by the DAC, upconverter, amplifier, and preceding transmission segments.  Therefore, the phase offset between the commanded RF drive and the forward drive may be filtered to directly 
infer the phase drift $\varphi^{\mathrm{fwd}}(t)$ accumulated along the forward path. This procedure isolates phase rotation introduced 
by the hardware drive chain from phase shifts arising from physical 
detuning within the cavity.

\subsection{Cavity Voltage and Detuning}

Detuning plays a critical role in RF cavity performance because even small frequency offsets from resonance can significantly impact power efficiency and system stability. Detuning relative to the reference is defined by
\[
\Delta\omega(t) = \omega_r(t) - \omega^{\mathrm{ref}},
\]
where $\omega_r(t)$ is the instantaneous cavity resonance frequency. When the cavity is on resonance, reflected power is minimized for the given coupling. As the magnitude of detuning increases, the cavity impedance becomes increasingly reactive, resulting in higher reflected power and greater forward RF power required to maintain a specified cavity voltage. Consequently, sustained operation off resonance can quickly compromise high-power amplifier limits. At LANSCE, measurements of reflected power are routed to the resonance control (RC) system and the RF interface control (RFIC) chassis for monitoring. The RC system mechanically drives tuning slugs to minimize reflected power and maintain operation near resonance.

In the rotating frame, the complex envelope is accurately modeled by the first-order system given by \cite{kwon2011fpga,simrock2022low}
\begin{equation}
    \dot V^{\mathrm{cav}}
=
\left(-\omega_{1/2} + j\Delta\omega\right) V^{\mathrm{cav}}
+ \omega_{1/2} V^{\mathrm{fwd}}
+ d,
\end{equation}
where $\omega_{1/2} = \omega_0/2Q_L$ is the loaded
half-bandwidth and $d(t)$ is a complex disturbance normalized consistently with the
voltage definition.  This disturbance represents
the aggregate effects of beam loading \cite{scheinker2015adaptive}, calibration error, and residual
mismatch.  The reflected waveform is digitized 
and digitally downconverted in the FCM to produce complex IQ measurements
consistent with the cavity pickup and forward signals.  These measurements may be used to establish a systematic method to estimate detuning and distinguish it from hardware phase drifts accumulated along the drive and receiver chains. Incorporating the estimate into the control 
architecture not only improves drive prediction but increases operator awareness of 
developing faults that manifest as increasing reflected power.

\subsection{Pickup and Receiver Chain}

The cavity field is monitored by a weakly coupled pickup probe inside the cavity with which a signal proportional to the internal accelerating field is measured. The pickup analog RF signal is downconverted to IF near the cavity input coupler to reduce transmission loss and phase drift prior to transport to the FCM. The IF analog signal is converted by an ADC (analog-to-digital converter) to digital IF, filtered, and digitally downconverted to baseband in the FCM to produce complex IQ phasor signals that signify the amplitude and phase of the cavity field relative to the master reference. These IQ signals are supplied to the FPGA to compute the feedback drive command that regulates the cavity field and closes the control loop.

Any slow electrical length change or phase drift in the analog receiver chain resulting from mixers, cables, filters, and ADC drifts appear as a complex rotation of the measured baseband signal. In the complex envelope representation, the receiver drift is described by a phase shift $\varphi^{\mathrm{rec}}(t)$ between the true cavity field $V^{\mathrm{cav}}(t)$ and the measured pickup phasor $y^{\mathrm{cav}}(t)$. Because this rotation occurs downstream of the cavity and upstream of digitization with no additional phase references across the LLRF control architecture, the phase shift cannot be inferred from receiver data alone. The pickup measurement is modeled as
\[
y^{\mathrm{cav}}(t) = e^{j\varphi^{\mathrm{rec}}(t)} V^{\mathrm{cav}}(t) + \nu^{\mathrm{cav}}(t),
\]
where $\nu^{\mathrm{cav}}(t)$ represents noise from quantization and calibration errors. Consequently, $\varphi^{\mathrm{rec}}(t)$ must be treated as an unmeasured parameter or additional state variable that requires observer support.

\subsection{State-Space Formulation}

The IQ planar rotation matrix of a sufficiently small angle $\theta(t)$ is approximated by
\begin{equation}
    R(\theta)=  I+\theta J, \quad I=\begin{bmatrix}1&0\\0&1\end{bmatrix}, \quad J=\begin{bmatrix}0&-1\\1&0\end{bmatrix}.
\end{equation}
The state-space representation of the complex envelope is expressed in terms of the real vectors $x=[\Re V^{\mathrm{cav}}, \Im V^{\mathrm{cav}}]'$ and  $u=[\Re V, \Im V]'$ according to
\begin{eqnarray}
    \dot{x}(t) &=& A(t)
x(t)
+
B(t)u(t)
+
d(t), \label{eq:cont_baseband} \\
y(t) &=& C(t) x(t) + v(t)
\end{eqnarray}
where state, control, and measurement matrices are given by
\begin{eqnarray}
A = \Delta\omega J-\omega_{1/2} I, \; B = \omega_{1/2}R(\varphi^{\mathrm{fwd}}), \; C = R(\varphi^{\mathrm{rec}}).
\end{eqnarray}
Digital implementation is obtained by integrating the above model with a forward Euler step and uniform sampling period $T_s$.  Over the discrete time interval $[k,k+1)$, the evolution of the cavity baseband state is approximated by
\begin{eqnarray}
x_{k+1} &=& A_k x_k + B_k u_k+ d_k, \\
y_{k} &=& C_{k} x_{k} + \nu_{k},
\end{eqnarray}
where indices represent the discrete-time labels.  In operation at LANSCE, the commanded input $u_k$ is held constant over each time segment $[k,k+1)$ of the feedback sequence.  Although matrix exponentials may be derived in closed form to represent the exact evolution under piecewise constant variables, we assume that the sampling frequency is sufficiently fast to accurately model dynamics with first-order sampling methods.  Besides, any residual mismatch is lumped into the additive disturbance. 

\section{Informed Observation}
\label{sec:estimation}

Section~\ref{sec:baseband-model} expresses the drive chain, cavity, and pickup
channels in a common baseband frame. The dominant phase uncertainties are well represented by drift in the forward chain that rotates the commanded
actuation before it reaches the cavity coupler, cavity
resonance during the macropulse, and
 drift in the receiver chain that rotates the true
cavity field before it appears in the pickup channel.  This section establishes a framework to estimate the evolution of these phase uncertainties by integrating existing measurements with informed models.  Assume that a drive signal has been applied to the cavity over the interval
$[k-1,k)$ and that a new cavity field measurement becomes available at
time $k$. The observers in this section create a prediction using quantities
available before processing the new measurement, forms an innovation, and
finally updates the internal estimates for the next propagation interval
$[k,k+1)$.

\subsection{State and Disturbance Observer}
\label{subsec:unified-eso}

A reliable method for robust regulation is to
estimate a small set of effective unknowns than to attempt full physical
identification \cite{qiu2015application,qiu2021application}. For comparison to our proposed observer, we first consider a standard extended state observer (ESO) that
tracks the cavity field while adapting an additive disturbance that lumps uncertainties appearing within the baseband model \cite{vincent2011active}.  Let
$\hat A_{k-1}$, $\hat B_{k-1}$ denote the state and control matrices used for prediction over
$[k-1,k)$. The predicted state is given by
\begin{equation}
\label{eq:predicted_state_unified}
\hat x_k^{\mathrm{pred}}
  = \hat A_{k-1}\hat x_{k-1}
  + \hat B_{k-1}u_{k-1}
  + \hat d_{k-1},
\end{equation}
and the innovation is
\begin{equation}
\label{eq:innovation}
r_k = y_k - \hat x_k^{\mathrm{pred}}.
\end{equation}
The ESO update is
\begin{equation}
\label{eq:eso-base-update}
\hat x_k = \hat x_k^{\mathrm{pred}} + \alpha_x\, r_k,
\qquad
\hat d_k = \hat d_{k-1} + \alpha_d\, r_k,
\end{equation}
where $\alpha_x\ge 0$ and $\alpha_d\ge 0$ are gains that set the adaptation bandwidth. In
practice, these gains are selected to ensure that the linearized ESO error
dynamics are stable, robust to
measurement noise, and fast compared to slow drifts. Equations~\eqref{eq:predicted_state_unified}-
\eqref{eq:eso-base-update} define the ESO or standard observer.

\subsection{Forward Drift Observer}
\label{subsec:fwd-drift-observer}

Forward drift is inferred from the difference between the planar
orientations of the measured forward drive $u^{\mathrm{fwd}}_k$ and the
commanded drive $u_k$.
The raw forward drift measurement is defined by
\begin{equation}
\tilde{\varphi}^{\mathrm{fwd}}_k
=
\angle u^{\mathrm{fwd}}_k
-
\angle u_k,
\end{equation}
where $\angle z = \tan^{-1}(z_y/z_x)$ denotes the planar orientation of a vector $z$.  Because hardware drift evolves slowly relative to the control bandwidth, the estimate is updated with the low-bandwidth integrator defined by
\begin{equation}
\hat{\varphi}^{\mathrm{fwd}}_k
=
\hat{\varphi}^{\mathrm{fwd}}_{k-1}+\alpha_{\mathrm{fwd}}
\tilde{\varphi}^{\mathrm{fwd}}_k.
\end{equation}
where $\alpha_{\mathrm{fwd}}\ge 0$ is the observer gain.
The estimated forward drift is applied to the
drive command through the control matrix update so that the cavity input remains aligned with the
controller reference despite hardware phase rotation.

\subsection{Detuning Observer}
\label{subsec:detuning-observer}

The resonance frequency of the cavity perturbs with wall deformation caused by Lorentz forces, thermal expansion of the cavity structure, and microphonic excitation from vibrations in the cryomodule and surrounding vacuum vessel \cite{simrock2005digital}. In practice, such detuning often contains multiple vibration modes, bias offsets, and slow frequency wander that are difficult to represent with a fixed mechanical model.  For inference during closed loop regulation, we directly estimate detuning from the reflected wave measurement using a reflected wave model infused with a low bandwidth integrator.
At the input coupler, the measured reflected wave $u^{\mathrm{refl}}$ is the direct reflection of the incident forward wave and radiation from the stored cavity field leaking back through the same port according to the approximation $u^{\mathrm{refl}} \approx u^{\mathrm{fwd}} - \kappa x$,
where $\kappa$ describes the coupling of the internal cavity field into the reflected measurement channel.

In steady operation, by solving for $x$ in \eqref{eq:cont_baseband}, neglecting additive disturbance for readability, and substituting the resulting relation into the linear port model 
 gives the reflected phasor estimation
 \begin{eqnarray}
\label{eq:port_proxy_int}
\hat u^{\mathrm{refl}}_k
=
 u^{\mathrm{fwd}}_{k-1}
-\kappa \left(I- \frac{\hat \Delta \omega_{k-1}}{\omega_{1/2}} J\right)^{-1}u^{\mathrm{fwd}}_{k-1}.
\end{eqnarray}
Detuning is estimated according to the value with which the reflected wave measurement $u^{\mathrm{refl}}_k$ and the model estimate $\hat u^{\mathrm{refl}}_k$ agree.  We use gradient decent on the squared norm of the reflected innovation $r^{\mathrm{refl}}_k = u^{\mathrm{refl}}_k - \hat u^{\mathrm{refl}}_k$ to search for a detuning estimate that minimizes the norm of the innovation according to
\begin{equation}\label{eq:gd_update}
\hat\Delta\omega_k
= \hat\Delta\omega_{k-1}
+ \alpha_{\omega}
(r^{\mathrm{refl}}_k)' J u^{\mathrm{fwd}}_{k-1},
\end{equation}
where $\alpha_{\omega}\ge 0$ is the observer gain or step size normalized to the gradient of the objective.  Here, we have approximated the gradient to first order with respect to detuning.  Note that this update defines a smooth integrator tailored to the sensitivity of cavity resonance with respect to the reflected wave.
The detuning estimate is used immediately in the subsequent cavity
prediction by updating $\hat A_k=A(\hat\Delta\omega_k)$ for the next step in \eqref{eq:predicted_state_unified}.

\subsection{Receiver Drift Observer}
\label{subsec:rec-drift-observer}

Receiver drift represents a slow phase
rotation applied to the cavity pickup phasor after the cavity dynamics and
before digitization. Because this rotation cannot be inferred from monitor comparisons, it is
estimated by comparing the measured pickup phasor
$y_k$ with the baseband model prediction $\hat x^{\mathrm{pred}}_k$. The raw receiver drift measurement is obtained from the difference
between the orientations of the measured and predicted phasors according to
\begin{equation}
\tilde{\varphi}^{\mathrm{rec}}_k
=
\angle y_k
-
\angle \hat x^{\mathrm{pred}}_k.
\end{equation}
As with the forward drive phase drift,
the estimate of the receiver drift is updated using the low-bandwidth filter defined by
\begin{equation}
\hat{\varphi}^{\mathrm{rec}}_k
=
\hat{\varphi}^{\mathrm{rec}}_{k-1}+\alpha_{\mathrm{rec}}
\tilde{\varphi}^{\mathrm{rec}}_k,
\end{equation}
where $\alpha_{\mathrm{rec}}\ge0$ is the observer gain.  When the receiver drift observer is activated and used in feedback control, the state and additive disturbance observers in \eqref{eq:eso-base-update} are computed following the computation of $\hat{\varphi}^{\mathrm{rec}}_k$.  In particular, the cavity pickup and state estimate used to define the innovation in \eqref{eq:innovation} are first aligned to the same reference by either inversely rotating the pickup signal to the controller reference frame or by rotating the state estimate to the drifted reference frame. 

\section{Control Architecture}
\label{sec:control-architecture}

The objective of the LLRF controller is to regulate the cavity field to track a complex envelope $x^*_k$ while 
compensating for effects of detuning, beam loading, slow phase drift, and measurement noise.  Compensation is based on feedforward and feedback control models in which operational artifacts and disturbances are estimated online and used in the underlying models for improved tracking performance.  In this section, we introduce the feedforward and feedback methods under the assumption that estimates of model parameters are provided by the observers defined in the previous section.

\subsection{Feedforward}
\label{subsec:feedforward}

To drive the state estimation to the setpoint in one step, substitute $\hat x_{k+1}^{\mathrm{pred}} = x_{k+1}^*$ into \eqref{eq:predicted_state_unified} and invert $\hat B_k$ to obtain the feedforward command \cite{kwon2022iterative}
\begin{equation}
u_k^{\mathrm{ff}}
  = \hat B_k^{-1}
    \left( x_{k+1}^* - \hat A_k \hat x_k -  \hat d_k \right).
\label{eq:ff-law}
\end{equation}
The control matrix is invertible for all practical (nonzero) values of $\omega_{1/2}$.  If the estimated model agrees exactly with the plant, then this control command achieves perfect tracking in the absence of measurement and process noise.

\subsection{Feedback}
\label{subsec:feedback}

In addition to the feedforward component, the
controller constructs the commanded input
\begin{equation}
 u_k =  u_k^{\mathrm{ff}} + u_k^{\mathrm{fb}},
\qquad
 u_k^{\mathrm{fb}} = K_k \left(x_k^*-\hat x_k \right),
\label{eq:fb-law}
\end{equation}
where $K_k$ is a feedback gain.  The feedback term in
\eqref{eq:fb-law} is derived by linear quadratic regulation over one time step.  The objective is to attenuate deviations from the desired
regulation target while moderating control effort. Consider the quadratic performance index defined by the expression
\begin{equation}
\hat e_{k+1}' Q\, \hat e_{k+1} +  u_k' R\, u_k,
\label{eq:onestep-cost}
\end{equation}
where $Q$ and $R$ are positive semidefinite matrices, with one strictly positive definite, and $\hat e_k=\hat x_k - x_k^*$ is the estimated error.  Substituting the error dynamics
into \eqref{eq:onestep-cost} and minimizing with respect to $ u_k$ leads to the optimal gain matrix
\begin{equation}
K_k = \left(\hat B_k' Q \hat B_k + R\right)^{-1} \hat B_k' Q \hat A_k.
\label{eq:onestep-K}
\end{equation}
This gain formula may be expressed in closed form using properties of 2 by 2 matrices to enable relatively simple implementation onto FPGA firmware.

\section{Regulation and Fault Localization}
\label{sec:mc_fault_awareness}

A Monte Carlo study is performed to evaluate observer robustness and the ability of the proposed estimation architecture to correctly attribute error to its physical source. In RF cavity systems, multiple mechanisms can produce similar phase or amplitude deviations in the baseband response.  A controller that regulates the field well but misidentifies the origin of these perturbations provides limited diagnostic value. The MC analysis therefore evaluates both regulation performance and the reliability with which the observer assigns deviations to the correct subsystem. Two estimation architectures are compared. The proposed observer explicitly estimates forward drift, receiver drift, detuning, and a residual disturbance within the baseband model according to the observers defined above. The standard observer defined solely by equations \eqref{eq:eso-base-update} is also implemented for comparison.

\subsection{Monte Carlo Construction}
\label{subsec:mc_construction}

\begin{figure}
    \centering
    \includegraphics[width=\linewidth]{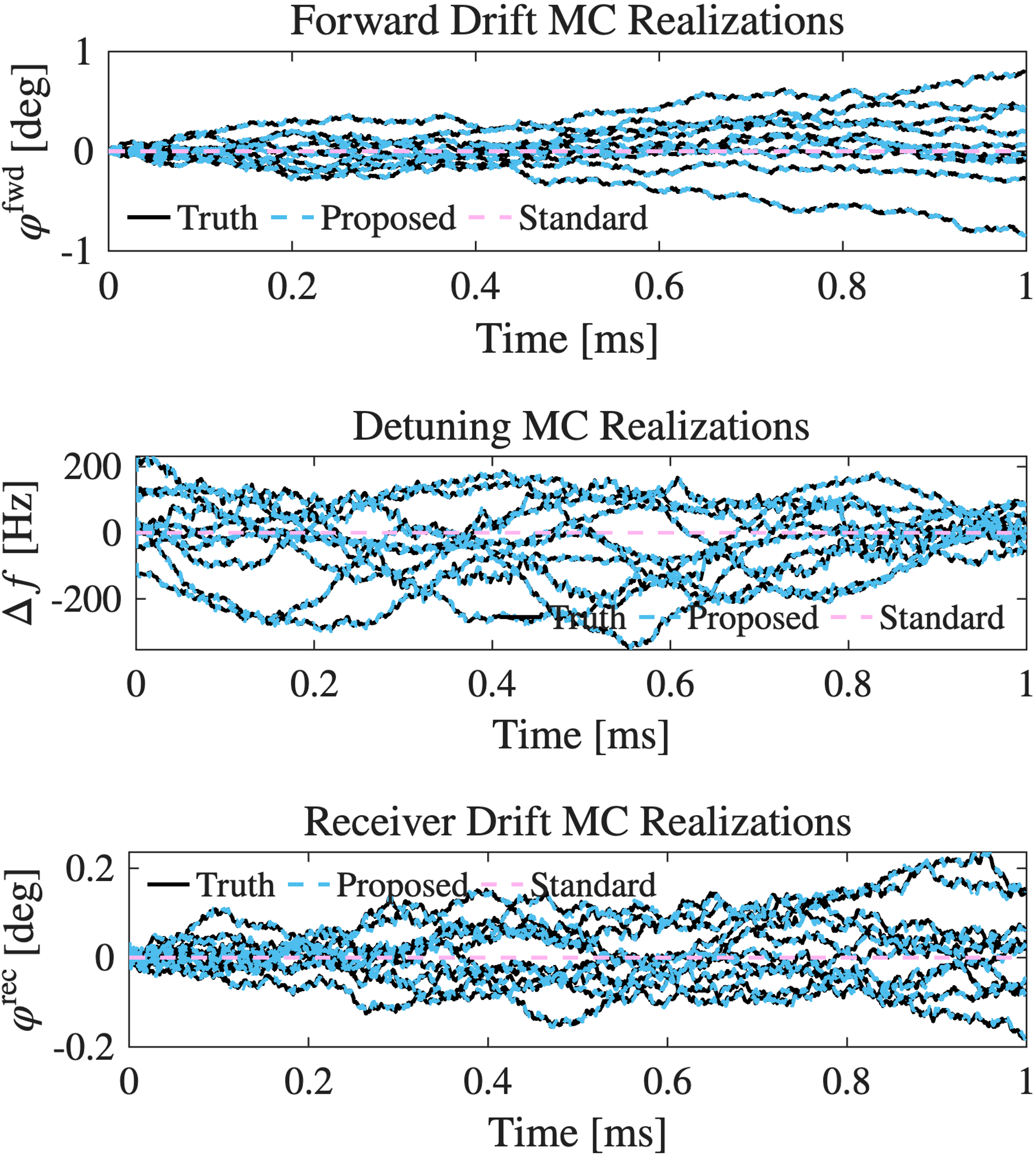}
    \caption{Representative MC realizations (showing 10/10k) of true and estimated forward-path phase drift, detuning, and receiver-path phase drift for the proposed and standard observers. Realizations are depicted over the fill and flattop portions of the macropulse.}
    \label{fig:mc_drifts}
\end{figure}

Each MC realization simulates the discrete cavity model and control architecture over the fill and flattop of a macropulse spanning 1 ms.  The simulations use ten thousand MC realizations to ensure sufficient exposure of systematic differences between the observer architectures.  System parameters are documented in Tables \ref{tab:system_param} and \ref{tab:observer_param}. These are fixed across all MC trials for both the proposed and standard observers.  For each MC trial, new randomized truth profiles are generated for the drift channels and noise processes. Detuning is generated using a composite process consisting of a bias term, several sinusoidal components with randomized frequency and amplitude, low-frequency stochastic wander, and a thermal contribution derived from mixer temperature variations. Phase drifts in the forward and receiver chains are generated by slowly evolving stochastic trajectories with small periodic components to represent variations in the RF drive and receiver chains.  Trajectories of the additive disturbance vector in the baseband model are generated by similar random processes. Gaussian measurement noise is injected into the cavity pickup, reflected-wave measurement, reference, and forward-drive channels each with a standard deviation of 1e-4.  Ten of the ten thousand MC drifts and detuning trajectories as well as their estimations using both the proposed and standard observers are shown in Figure \ref{fig:mc_drifts}.

\begin{table}
\centering
\begin{tabular}{c c c c c}
\hline
$T_s$ & $Q_L$ & $\omega_0$ & $Q$ & $R$ \\
\hline
1e-6 & $1.61$e4 & $2\pi \times 805$e6 & $0.1I$ & $100I$ \\
\hline
\end{tabular}
\caption{Simulation parameters for regulating the baseband voltage command of an 805 MHz side-coupled cavity.} \label{tab:system_param}
\end{table}

\begin{table}
\centering
\begin{tabular}{c c c c c c}
\hline
$\alpha_x$ & $\alpha_d$ & $\alpha_\omega$ & $\alpha_{\mathrm{fwd}}$ & $\alpha_{\mathrm{rec}}$  & $\kappa$ \\
\hline
$0.1$ & $1$e-4 & $1.0$ & $0.9754$ & $0.7316$ & $200$ \\
\hline
\end{tabular}
\caption{Parameters for propagating the proposed observer.  The standard observer uses $\alpha_x=0.1$ and $\alpha_d=0.3$, with remaining parameters set to zero.} \label{tab:observer_param}
\end{table}

\subsection{Regulation}
\label{subsec:mc_metrics}

\begin{figure}
    \centering
    \includegraphics[width=\linewidth]{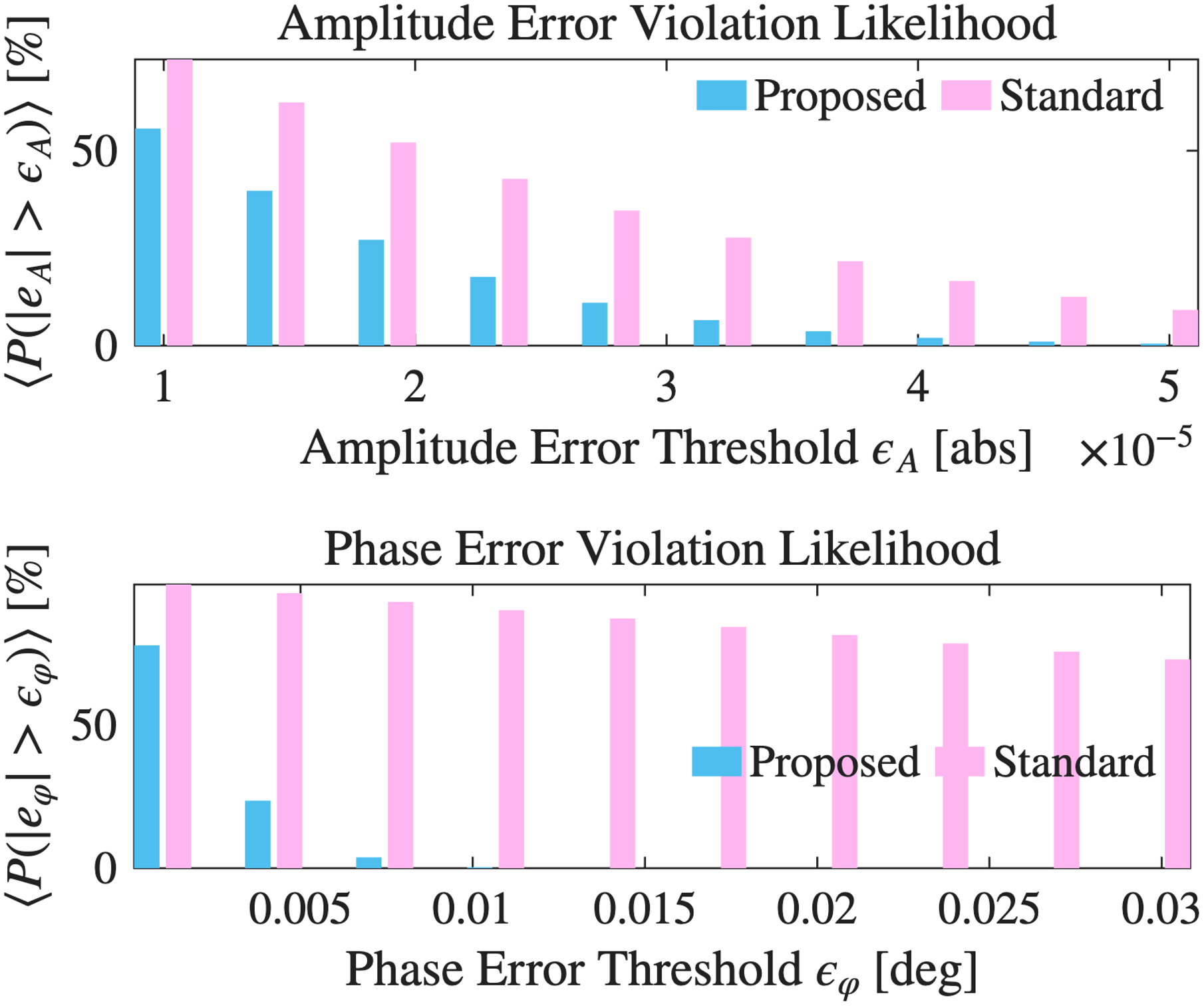}
    \caption{Likelihood that the average amplitude and phase regulation errors over time exceed the thresholds specified.}
    \label{fig:state_likelihood}
\end{figure}

Cavity field regulation is quantified using amplitude and phase errors between the cavity phasor and the reference trajectory.  The amplitude and phase of the cavity field vector are defined by $A(t) = \sqrt{x_1^2(t)+x_2^2(t)}$ and $\varphi(t)=\text{tan}^{-1}(x_2/x_1)$, respectively.  Amplitude and phase errors between the cavity and reference signals are defined by
\begin{equation}
    e_A(t) = A(t) - A^*(t), \qquad e_{\varphi}(t) = \varphi(t)-\varphi^*(t),
\end{equation}
where the asterisk denotes the reference signature.  Figure \ref{fig:state_likelihood} displays the probability of the MC trials averaged over time to violate amplitude and phase error margins by the specified thresholds.  Both the proposed and standard observers demonstrate the ability to regulate the cavity field while rejecting the effects of disturbances and noise.  One distinct outcome of the proposed observer is the ability to maintain tight phase margins in the presence of receiver phase drift.  In contrast, the standard observer does not compensate for errors caused by the receiver chain because its extended state only appears as an additive term in the baseband model.  Let us note that further simulations show that phase errors incurred by the standard observer tend to decrease to the margins of the proposed observer when disturbance drifts in the receiver channel are bypassed by setting $\varphi^{\mathrm{rec}}(t)=0$.

\subsection{Localized Fault Identification}
\label{subsec:mc_results}

\begin{figure}
    \centering
    \includegraphics[width=\linewidth]{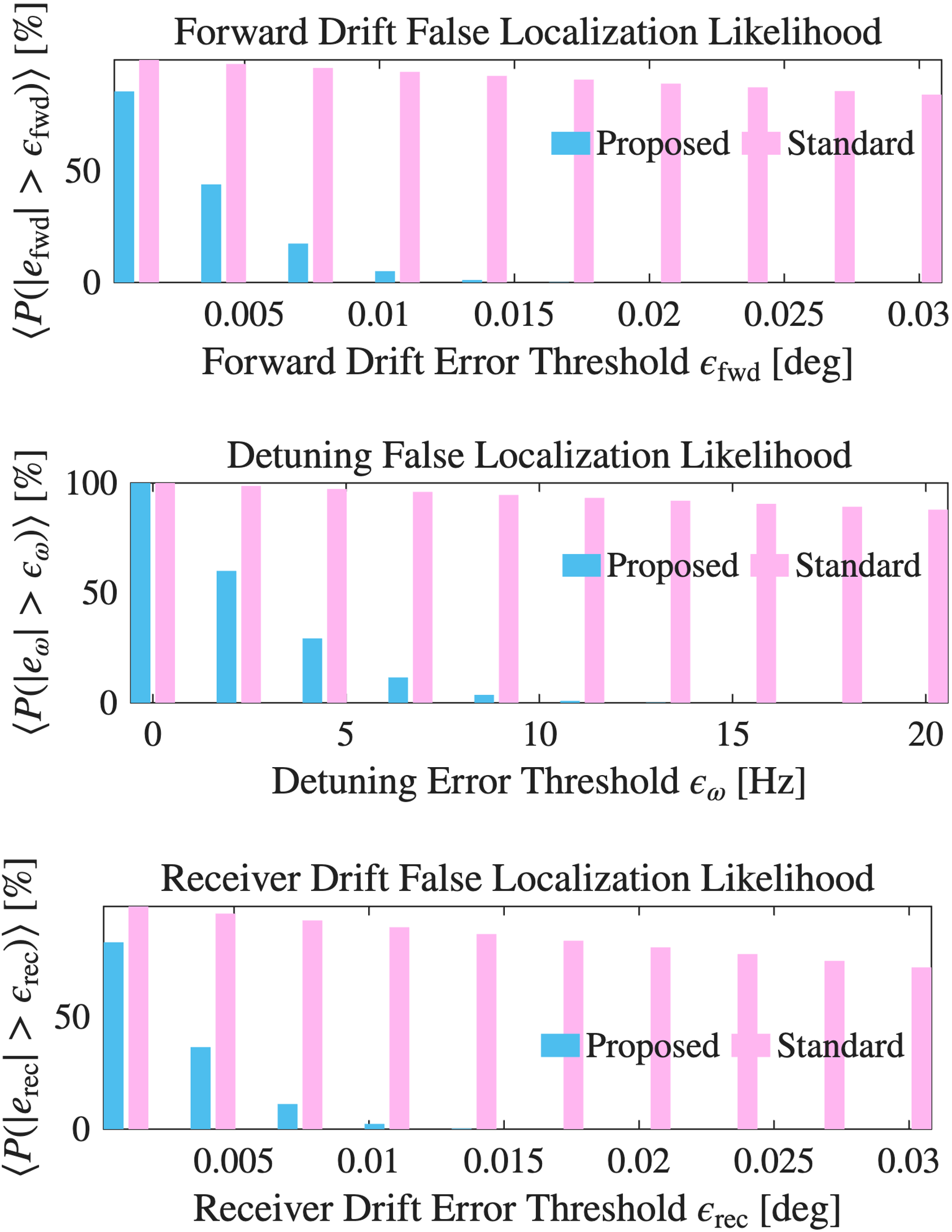}
    \caption{Likelihood that the average error over time between the true and estimated disturbances exceeds the specified thresholds.}
    \label{fig:drift_likelihood}
\end{figure}

Awareness of faults is evaluated by measuring the ability of the observer to correctly assign deviations to their physical source. A false-localization event is defined when the estimated drift of a given channel differs from the corresponding true drift by more than a specified threshold. For each channel, the fraction of MC realizations exceeding a given threshold averaged over time defines the likelihood of false localization. In analogy with the amplitude and phase violation metrics shown in Figure \ref{fig:state_likelihood}, the same probability construction is applied here to estimation error relative to the true drift trajectories. The resulting likelihoods are summarized in Figure \ref{fig:drift_likelihood} for detuning and drifts in the forward and receiver chains. 

Figure \ref{fig:drift_likelihood} demonstrates a clear separation between the two observer architectures in their ability to localize faults. The proposed observer maintains low false-localization likelihood across all channels. This indicates that deviations in the cavity response are consistently attributed to the correct subsystem even under stochastic drift and measurement noise.  Because forward drift, receiver drift, and detuning enter through distinct measurement relationships in the proposed observer formulation, their effects remain separable in the estimated states as the underlying trajectories evolve. In contrast, the standard observer frequently exhibits elevated false-localization likelihood because the effects across channels are absorbed into a single disturbance estimate. Consequently, physically distinct mechanisms such as drive-chain phase drift, receiver-path variation, and cavity detuning are not reliably distinguished. Although both observer architectures maintain acceptable cavity field error margins under diverse operating conditions, the diagnostic interpretability of their estimated states differs significantly.

\section{Conclusion}
\label{sec:conclusion}

This work has presented an observer framework for RF cavity field regulation that explicitly separates several disturbance mechanisms commonly encountered in accelerator RF systems. By augmenting a standard ESO with dedicated estimators for forward-chain phase drift, receiver-chain phase drift, and cavity detuning, the proposed architecture enables simultaneous regulation of the cavity field and identification of the physical origin of phase perturbations. Monte Carlo analysis demonstrates that the proposed observer achieves regulation performance comparable to or better than conventional disturbance-based observers while significantly improving the reliability with which drift mechanisms are identified. In particular, simulations show an exponential reduction in false localization of detuning, forward-drive phase drift, and receiver-chain drift when compared with a standard ESO observer that models these effects as a single lumped disturbance. This improvement arises because each disturbance channel is inferred through a measurement relationship that is physically consistent with the underlying RF system.

The ability to localize disturbances has important implications for accelerator operations. In many facilities, increases in reflected power or phase error can arise from multiple sources that are difficult to distinguish in real time. By providing separate estimates for cavity detuning and hardware phase drift, the observer enables operators and diagnostic systems to identify the most likely subsystem responsible for a deviation. Such information can guide corrective actions such as resonance control adjustments, RF drive chain inspection, or calibration of receiver electronics. In large accelerator facilities where downtime and recovery procedures can be costly, improved diagnostic clarity can significantly reduce the time required to identify and resolve RF system faults.

Several directions for future work remain. Experimental validation of the proposed observer on operational RF stations would provide valuable insight into the interaction between estimator dynamics and real hardware nonlinearities. Additional extensions could consider mechanical models of cavity detuning, including microphonic vibration modes and Lorentz-force effects, which are particularly relevant in superconducting accelerators \cite{padamsee2014superconducting,liepe2001dynamic,simrock2005digital,yamamoto2009observation}.  Stability of the feedback control should be rigorously established, and adaptive observer gains may be considered for operations that require even tighter stability margins. Finally, integration of the observer outputs into higher-level machine protection and supervisory control systems may enable automated detection of developing RF faults before they impact accelerator performance.  As accelerator facilities continue to operate at higher beam power and tighter field stability requirements, the ability to regulate RF cavities while simultaneously diagnosing the origin of disturbances will become increasingly important. Observer-based approaches that explicitly model the underlying physics of RF systems provide a promising path toward achieving both goals.

\linespread{1}
\bibliographystyle{unsrt}
\bibliography{references.bib}

\end{document}